# Monte Carlo Dynamically Weighted Importance Sampling For Finite Element Model Updating

**Daniel J Joubert, Tshilidzi Marwala.**
University of Johannesburg, Kingsway Campus, Cnr. Kingsway and University Rd, Auckland Park

**ABSTRACT**

The Finite Element Method (FEM) is generally unable to accurately predict natural frequencies and mode shapes of structures (eigenvalues and eigenvectors). Engineers develop numerical methods and a variety of techniques to compensate for this misalignment of modal properties, between experimentally measured data and the computed result from the FEM of structures. In this paper we compare two indirect methods of updating namely, the Adaptive Metropolis Hastings and a newly applied algorithm called Monte Carlo Dynamically Weighted Importance Sampling (MCDWIS). The approximation of a posterior predictive distribution is based on Bayesian inference of continuous multivariate Gaussian probability density functions, defining the variability of physical properties affected by forced vibration. The motivation behind applying MCDWIS is in the complexity of computing normalizing constants in higher dimensional or multimodal systems. The MCDWIS accounts for this intractability by analytically computing importance sampling estimates at each time step of the algorithm. In addition, a dynamic weighting step with an Adaptive Pruned Enriched Population Control Scheme (APEPCS) allows for further control over weighted samples and population size. The performance of the MCDWIS simulation is graphically illustrated for all algorithm dependent parameters and show unbiased, stable sample estimates.

Introduction
Physical properties, for example, areas, moments of inertia, or elasticity (Young's Modulus), all influence the outcome of eigenvalues (natural frequencies) and the eigenvectors (mode shapes) of the structure. The method of model updating discussed in this paper is of the indirect or iterative type, thus taking into account changes of physical parameters when the structure is under forced vibration or dynamic load. This results in the mass and stiffness matrices of the updated FE model having physical meaning, i.e. eigenvalues and eigenvectors [1]. In recent years, it has been well established that there exists Bayesian evidence in FE model updating of structures, [1, 2 and 3]. Various MCMC sampling techniques have been applied [4], making predictions on uncertain or stochastic parameters in order to produce realistic modal frequencies of the models under vibration.

**1. Bayesian Inference and Gaussian Process Model**
It is important to note that due to the complexity of systems, a closed form posterior distribution is not analytically available and in most cases cannot be accurately approximated. This comes about through the multi-dimensionality of the search space in multimodal systems. This multi-dimensionality of the Monte Carlo integrals causes difficulties in computation, and leads to the use of approximation methods. Thus, for the purposes of this MCDWIS simulation the normalizing constant $Z(\boldsymbol{\theta})$ of the posterior distribution is not available and is referred to as analytically intractable. In Bayesian inference, the posterior is in fact generally considered to be analytically intractable [5]. The Bayesian approach is governed by Bayes Rule[6].

$$p(\boldsymbol{\theta}|D) = \frac{f(D|\boldsymbol{\theta})\pi(\boldsymbol{\theta})}{\pi(D)} \text{ or simply } p(\boldsymbol{\theta}|D) \propto f(D|\boldsymbol{\theta})\pi(\boldsymbol{\theta}) \tag{1}$$

Where $\boldsymbol{\theta}$ represents the vector of updating parameters, $f(D|\boldsymbol{\theta})$ is the likelihood probability distribution function and $\pi(\boldsymbol{\theta})$, denotes the prior probability distribution function. As shown, the posterior probability can be approximated through Bayesian inference. Since we will be considering a multi-dimensional parameter search space with non-independent variables in this study, we will be using continuous multivariate probability density functions. Gaussian Processes are very popular in Bayesian estimation procedures. And in this thesis, all probability functions will be Normal. The covariance of the target distribution is unknown. These statistical parameters are approximated analytically during the burn-in stages of the simulations. However, $\mu_0$ and $\Sigma_0$ are assumed to initialize the Gaussian process. These initializations are also known as hyper-parameters. The likelihood distribution is given as,

$$f(D|\boldsymbol{\theta}) = \frac{1}{((2\pi)^{N_m}|\Sigma_f|)^{0.5}} \exp\left(-0.5(f_{N_m}^c - f_{N_m}^m)^T \Sigma_f^{-1}(f_{N_m}^c - f_{N_m}^m)\right). \tag{2}$$

Here $N_m$ denotes the dimension of the frequency vector, $f_i^m$ and $f_i^c$ represent the measured and computed natural frequencies respectively. The determinant of the frequency covariance matrix is $|\Sigma_f|$. The Taylor expansion can be used to express the log-likelihood distribution function as,

$$\log(f(D|\boldsymbol{\theta})) = \frac{N_m}{2}\ln(2\pi) + 0.5\ln|\Sigma_f| + 0.5(f_{N_m}^c - f_{N_m}^m)^T \Sigma_f^{-1}(f_{N_m}^c - f_{N_m}^m). \tag{3}$$

The prior probability distribution is also chosen as a Gaussian probability density function.

$$\pi(\boldsymbol{\theta}) = \frac{1}{((2\pi)^d|\Sigma_{\theta,0}|)^{0.5}} \exp\left(-0.5(\boldsymbol{\theta} - \boldsymbol{\mu_0})^T \Sigma_{\theta,0}^{-1}(\boldsymbol{\theta} - \boldsymbol{\mu_0})\right) \tag{4}$$

Where $d$ denotes the dimension of the variable parameter vector, $\boldsymbol{\theta}$. The process is initialized by the parameter covariance matrix, $\Sigma_{\theta,0}$ and mean vector, $\boldsymbol{\mu_0}$. Consequently, given the likelihood and prior, it is obvious to note that multiplying two Gaussian distributions yields a Gaussian distribution. Therefore these distributions are conjugate distributions. From Bayes' theorem the posterior probability distribution can be approximated as,

$$p(\boldsymbol{\theta}|D) = \frac{1}{Z(\boldsymbol{\theta},D)} \exp\left(-0.5(f_{N_m}^c - f_{N_m}^m)^T \Sigma_f^{-1}(f_{N_m}^c - f_{N_m}^m) - 0.5(\boldsymbol{\theta} - \boldsymbol{\mu_{\theta n}})^T \Sigma_{\theta,n}^{-1}(\boldsymbol{\theta} - \boldsymbol{\mu_{\theta n}})\right). \tag{5}$$

Where the normalizing constant is described by the integral,

$$Z(\boldsymbol{\theta}, D) = \int f(D|\boldsymbol{\theta})\pi(\boldsymbol{\theta})d\boldsymbol{\theta}. \tag{6}$$

Due to the multi-dimensionality of this integral, solving it is considered computationally complex. Attempts have been made to estimate reasonable normalizing constants in model updating[6]. However, this is not effective given increasingly complex systems. This is to say that there is no 'one size fits all' approach in assigning normalizing constants to the posterior distribution for all structural models with varied sensitivity parameters and material properties.

## 2. Monte Carlo Dynamically Weighted Importance Sampling (MCDWIS)

A key advantage of MCDWIS is that this algorithm does not require perfect sampling. Thus, MCDWIS is interesting to consider with regard to parameter estimation in model updating of FE models. Due to the computationally expensive effort of solving to FE model at every time step during the simulation, model updating is inherently computationally very intensive. Avoiding the need for exact sampling which in itself can be very expensive or in most cases not possible, makes MCDWIS a worthwhile technique to experiment with in model updating, adopting a stochastic Monte Carlo approximation approach. The use of the MCDWIS algorithm is thus validated in approximating unbiased estimates which is controlled within a desired reduced variance range from the originally drawn samples. Reduced variance is simply through the use of an importance sampling estimate. New weighted samples are then generated, referred to in literature as the new population. Through population control schemes the system becomes adaptive to the requirements of the outcome, as the posterior distribution is left invariant with respect to dynamic importance weights.

The prior distribution denoted as $\pi(\boldsymbol{\theta})$, represent the prior knowledge through observation of the previous system state. From the observed auxiliary data $y$ the likelihood function for the statistical model is, [7]

$$f(y|\boldsymbol{\theta}) = \frac{p(\boldsymbol{\theta},y)}{Z(\boldsymbol{\theta})} \tag{7}$$

Where $p(\boldsymbol{\theta}, y)$ denotes the conditional probability density function of the state and auxiliary samples. $Z(\boldsymbol{\theta})$ is the normalizing constant which is dependent on the value of $\boldsymbol{\theta}$. The posterior probability distribution for variable $\boldsymbol{\theta}$ can thus be describe as [7],

$$p(\boldsymbol{\theta}|y) = \frac{1}{Z(\boldsymbol{\theta})} p(\boldsymbol{\theta},y)\pi(\boldsymbol{\theta}) \tag{8}$$

This algorithm, similar to importance sampling uses important weights to prioritize acceptance and locate the optimal search space for sampling. Each state in the Markov chain of size $N$ is given by the joint density, $(\boldsymbol{\theta}, \boldsymbol{w}) = \{\theta_1, w_1; \ldots; \theta_N, w_N\}$. Every move step of the algorithm involves two actions, namely:

1. Dynamic weighting: Where each state is updated by a dynamic weighting transition step in order to compute a new population. [8]
2. Population control: Samples related to small weights with regard to updated results from the FE model are discarded while weighted samples with stronger relevance to the model objective function are duplicated in the new population. This induces biased samples but is counter balanced by assigning new weights to the samples for every new population.

Figure 5 below demonstrates the concept of the APEPCS in the dynamic weighting and population control stages of the algorithm [8]. The samples' weight classification residing outside of the weight control bounds $W_{up}$ and $W_{low}$ are scaled by the use of the control parameter $d_s$ in the enriching stage where $w_{t,i}$ is too large and the probability density parameter $q$ during the pruning stage ($w_{t,i} < W_{low}$), of the population control scheme. Weighting is controlled by the ratio of the upper and lower weight control bounds and can be described as, $\kappa = \frac{W_{up}}{W_{low}}$. This ratio controls the moving ability of the system and is called the freedom parameter [9]. APEPCS controls the population size within the range $[n_{min}, n_{max}]$. In this fashion the algorithm constructs new populations with each run, adjusting the size according to weight assignments. These weights influence the acceptance of elements in the Markov Chain into the new population.

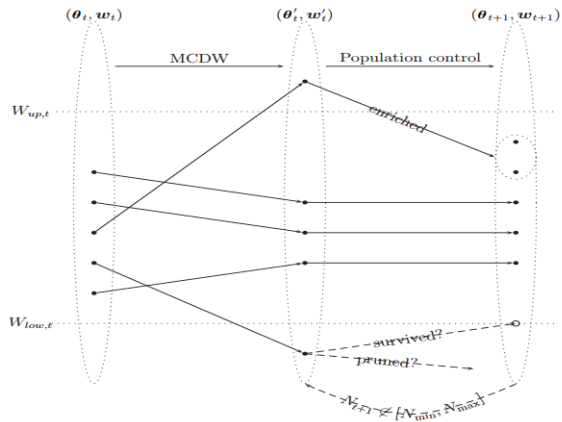

**Fig. 1 DWIS Algorithm**

The details on the different steps involved in the APEPCS process can be found in Liang[10]. There are two population control schemes of which for this problem the authors propose using R-scheme for simplicity of the algorithm. According to Liang [9], Scheme-R fits best with Monte Carlo simulations.

*Scheme –R*

1. Dynamic weighting: The scheme is applied to $(\theta_{t-1}, w_{t-1})$. Where $W_c$ denotes the dynamic move switching parameter of the value $\varphi_t$ between 0 and 1 depending on the value of, $W_{up,t-1}$. In the case where $W_{up,t-1} \leq W_c$ the value, $\varphi_t = 1$. Inversely when $W_{up,t-1} \geq W_c$, $\varphi_t$ is set to 0. $W_c$ is typically chosen as an exponential expression and relates to the weight control bounds by $W_c = 10^\kappa$. The new population is denoted as $(\boldsymbol{\theta_t'}, \boldsymbol{w_t'})$. The parameter $\varphi_t$ is chosen as a function of $(\boldsymbol{\theta_t}, \boldsymbol{w_t})$. It is interesting to note that, for the R-Type move, if $\varphi_t = 0$, the sampler is essentially a random walk process. And when $\varphi_t = 1$, then the sampler performs the R-Type move. Essentially, the purpose of introducing the $\varphi_t$ parameter is to prevent the weighting procedure from drifting to zero or infinity.
2. Population control: Next APEPCS is applied to the new population. And the new population becomes $(\boldsymbol{\theta_t}, \boldsymbol{w_t})$.

By applying this to model updating the aim is to see the MCDWIS algorithm produce more accurate results with faster convergence and see an improvement in computational efficiency. In order to evaluate the performance of the DW process we would like to see the behavior of the switching parameter, the adaptation of $W_{up,t}$, the adaptation of the population size and the log weight of selected $i^{th}$ states. To avoid unwanted inflation of the weighting process, it is important to implement a Weight Behavior Analysis (WBA). Further literature and proofs regarding weight behavior analysis can be found in Liang [10]. The main acceptance and rejection process can be summarized as follows:

1. Firstly a proposed sample $\theta^*$ is drawn from the proposal distribution $\pi(\theta^*|\theta_t)$.
2. The auxiliary simulated state trajectories from $f(y|\theta^*)$, are generated through MH a procedure $\boldsymbol{y} = \{y_1, \dots, y_M\}$. Next we compute the importance sampling (IS) estimate acting as the normalizing constant ratio $\hat{R}_t(\theta, \theta^*) = Z(\theta)/Z(\theta^*)$ [9]

$$\hat{R}_t(\theta_t, \theta^*) = \frac{1}{M}\sum_{j=1}^{M}\frac{p(y_j, \theta)}{p(y_j, \theta^*)} \tag{9}$$

3. From this the Monte Carlo dynamic weighting ratio is calculated. The conditional probability distribution can be expressed as, $p(\boldsymbol{y}, \theta_t) = p(\theta_t|\boldsymbol{y})\pi(\boldsymbol{y})$. [10]

$$r_d = r_d(\theta_t, \theta^*, w) = w\hat{R}_t(\theta_t, \theta^*)\frac{p(y,\theta^*)}{p(y,\theta_t)}\frac{\pi(\theta_t|\theta^*)}{\pi(\theta^*|\theta_t)} \tag{10}$$

4. As previously discussed, the joint density is updated using a uniformly drawn random numbers $U \sim unif(0,1)$ in the inequality, [10]. From APEPCS we determine $\varphi_t$

$$(\boldsymbol{\theta'}, \boldsymbol{w'}) = \begin{cases} \left(\boldsymbol{\theta^*}, \frac{r_d}{a}\right), & if\ U \leq a, \\ \left(\boldsymbol{\theta}, \frac{w}{1-a}\right), otherwise. \end{cases} \tag{11}$$

Where, $a = \frac{r_d}{r_d + \varphi_t}$.

Given the dynamic weights with associated sampled parameters, $(\theta_1, w_1), (\theta_2, w_2), \dots, (\theta_N, w_N)$ the weighted average of the samples can be expressed as [9]

$$\hat{\mu} = \sum_{t=burn+1}^{N}\sum_{i=1}^{n'}\left(\frac{w_{t,i}\rho(\theta_{t,i})}{w_{t,i}}\right) \tag{12}$$

Where $\rho(\theta_t)$ represents a state function over all the drawn samples. The mean value can be assumed to be consistent and asymptotically normally distributed. The proof for this can be found in Liang [10]. From simulating the structural model and updating the chosen sensitivity parameters behaving randomly in the system, graphs are plotted and conclusions are drawn to support the findings in performance. Various plots will test the validity of the simulation and their significance will be discussed where necessary in the coming chapters. The initial statistical parameters, constraints and initial conditions are decided upon and the various algorithm parameters are tuned by virtue of trial and error.

## 6. GARTUER SM-AG19 Simulation

The GARTEUR model is a well-known tool in the structural dynamics research field, providing researchers with common ground in testing different hypotheses in FE model updating. In this section, we simulate the MCDWIS algorithm and compare the results with an Adaptive Metropolis Hastings (AMH) algorithm. The model testbed being tested is specifically called the GARTUER SM-AG19 model, and more information can be found in papers by Degener [11] and Friswell [12]. The nominal length and width of the aero plane structure is 1.5 m and 3 m respectively. The aluminum structure has a total mass of 44 kg. A viscoelastic layer was bonded to the wings to induce an increased damping effect. The beam elements of the structure are all modelled as Euler-Bernoulli beam elements, with the materials assumed to be standard isotropic on all the elements. From Friswell [12], the natural frequencies with respect to the measured mode shapes are determined to be: 6.38, 16.10, 33.13, 33.53, 35.65, 48.38, 49.43, 55.08, 63.04, 66.52 Hz. The updating vector consists of the right wing moments of inertia and torsional stiffness $(RI_{min}, RI_{max}, RI_{torsional})$, the left wing moments of inertia and torsional stiffness $(LI_{min}, LI_{max}, LI_{torsional})$, the vertical tail moment of inertia $(VTP_{I_{min}})$ and the overall density of the structure $\rho$. The updating vector is given by, $\boldsymbol{\theta} = [\rho, VTP_{I_{min}}, LI_{min}, LI_{max}, RI_{min}, RI_{max}, LI_{torsional}, RI_{torsional}]$. The initial starting vector used is $\theta_0 = [2785, 8.34 \times 10^{-9}, 8.34 \times 10^{-9}, 8.34 \times 10^{-7}, 8.34 \times 10^{-9}, 8.34 \times 10^{-7}, 4 \times 10^{-8}, 4 \times 10^{-8}]$. With the maximum and the minimum bounds respectively used as, $\theta_{max} = [3500, 12 \times 10^{-9}, 11.2 \times 10^{-9}, 12 \times 10^{-7}, 11.2 \times 10^{-9}, 12 \times 10^{-7}, 6 \times 10^{-8}, 6 \times 10^{-8}]$ and $\theta_{min} = [2500, 6 \times 10^{-9}, 8 \times 10^{-9}, 6 \times 10^{-7}, 8 \times 10^{-9}, 6 \times 10^{-7}, 3 \times 10^{-8}, 3 \times 10^{-8}]$. With the diagonals of the covariance initialized as $\sigma = [5 \times 10^2, 5 \times 10^{-9}, 5 \times 10^{-9}, 5 \times 10^{-7}, 5 \times 10^{-9}, 5 \times 10^{-7}, 5 \times 10^{-8}, 5 \times 10^{-8}]$. Figure 2a to figure 2j demonstrates the mode shapes represented by the natural frequencies.

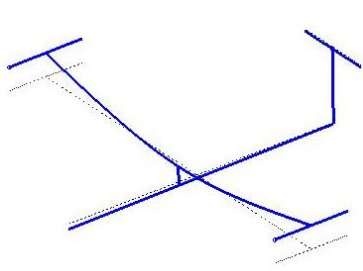 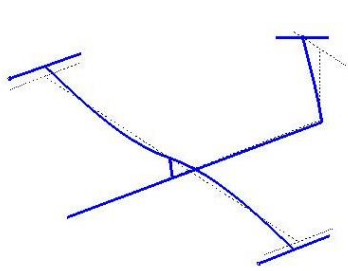 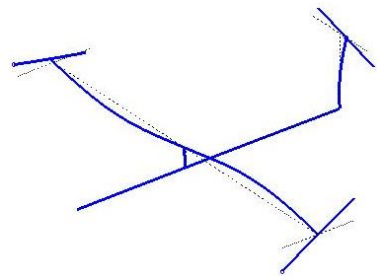

**Fig. 2a Mode 1 at 6.38 Hz**     **Fig. 2b Mode 2 at 16.10 Hz**     **Fig. 2c Mode 3 at 33.13 Hz**

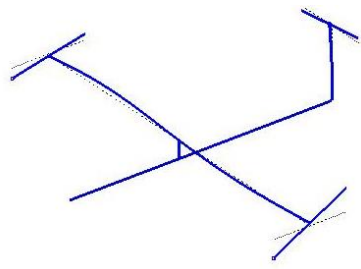 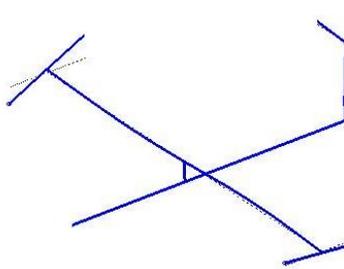 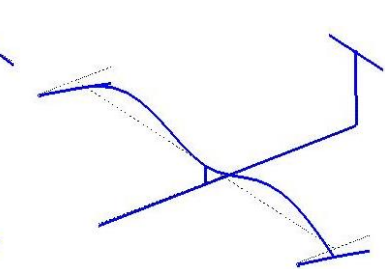

**Fig. 2d Mode 4 at 33.53 Hz**     **Fig. 2e Mode 5 at 35.65 Hz**     **Fig. 2f Mode 6 at 48.38 Hz**

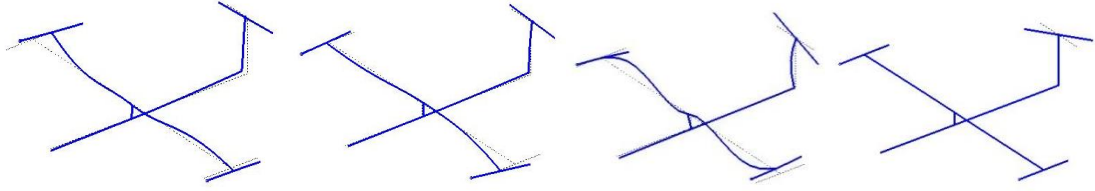

**Fig. 2g Mode7 at 49.43 Hz   Fig. 2h Mode 8 at 55.08 Hz   Fig. 2i Mode 9 at 63.04 Hz   Fig. 2j Mode 10 at 66.52 Hz**

## 6. Results

A fast run of the double MH algorithm is performed, which gives the initial set of samples for the MCDWIS to start off with. From 1000 iterations, the first 200 are discarded as burn-in and 200 samples are evenly selected from the remaining 800. This initial set of samples, are allocated with weight of 1 for initialization. The algorithm parameters are given in Table 1.

**Table 1: Algorithm Parameters**

| Algorithm parameter | Values |
|---|---|
| $N$ | 1000 |
| $burn-in$ | 250 |
| $W_c$ | $e^7$ |
| $n_{min}$ | 100 |
| $n_{max}$ | 1000 |
| $n_{low}$ | 200 |
| $n_{up}$ | 500 |
| $\lambda$ | 2 |
| $\kappa$ | $\log_{10} W_c$ |

The adaptation of the population size is given by Figure 3. Figures 4, shows a Gaussian representation of the log weight at the 10th state. Due to the higher dimensional vector, the algorithm reacts by allocating exceptionally high weights during the first few time steps and exponentially decreases, and then fluctuates. Figure 5 shows the adaptation of the upper bound weight value. Figure 6 shows the adaptation of the R-Scheme switching parameter $\varphi_t$, during computation. When the ratio between the highest and lowest weights exceed that of the allowable ratio between the upper and lower control bounds, the algorithm reacts accordingly, by gradually drifting back into the desired range, enabling the R-Scheme procedure to take effect. Figure 7 shows the correlation between the new samples at the 100th state, and the symmetry over the diagonal shows that the samples are correlated.

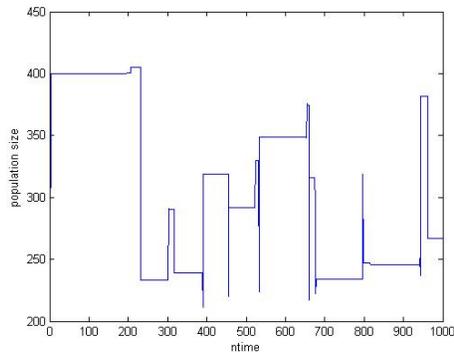
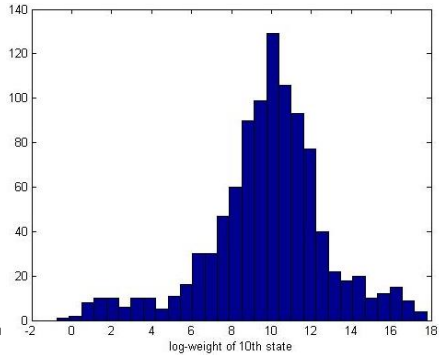

**Fig. 3 Population size adaptation**     **Fig. 4 Log weight at the 10th state**

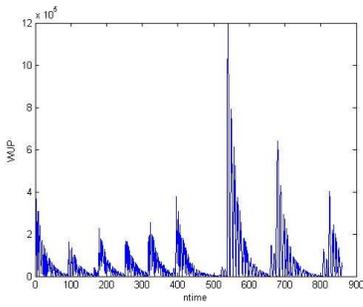 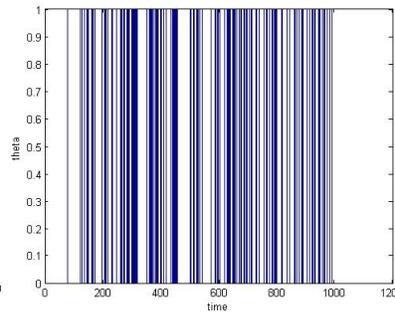 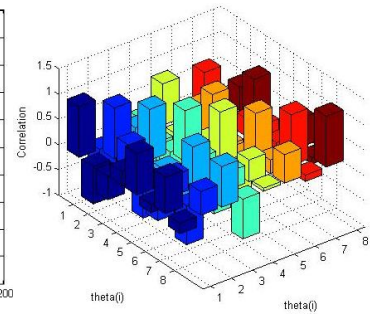

**Fig. 5 Adaptation of WUP**  **Fig. 6 Adaptation of the switching parameter**  **Fig. 7 Correlation**

Table 2 gives the values of the elements in the updated vector after the simulation. Table 3 gives the resultant natural frequencies corresponding to each respective mode shape.

**Table 2: Updating Vector results**

| d | Initial $\theta_0$ | $\theta$ vector, AMH method | $\theta$ vector, MCDWIS method |
|---|---|---|---|
| $\rho$ | 2785 | 2847 | 2736 |
| $VTP_{I_{min}}$ | $8.34 \times 10^{-9}$ | $9.1275 \times 10^{-9}$ | $7.1213 \times 10^{-9}$ |
| $LI_{min}$ | $8.34 \times 10^{-9}$ | $9.2278 \times 10^{-9}$ | $1.0200 \times 10^{-8}$ |
| $LI_{max}$ | $8.34 \times 10^{-7}$ | $7.6652 \times 10^{-7}$ | $7.9716 \times 10^{-7}$ |
| $RI_{min}$ | $8.34 \times 10^{-9}$ | $9.6244 \times 10^{-9}$ | $1.0202 \times 10^{-8}$ |
| $RI_{max}$ | $8.34 \times 10^{-7}$ | $8.3867 \times 10^{-7}$ | $6.0834 \times 10^{-7}$ |
| $LI_{torsional}$ | $4 \times 10^{-8}$ | $4.4388 \times 10^{-8}$ | $4.1102 \times 10^{-8}$ |
| $RI_{torsional}$ | $4 \times 10^{-8}$ | $3.5720 \times 10^{-8}$ | $3.6212 \times 10^{-8}$ |

**Table 3: Natural frequency results**

| Mode | Measured Frequency (Hz) | Initial Frequency (Hz) | Err (%) | AMH Frequency | Err (%) | MCDWIS Frequency | Err (%) |
|---|---|---|---|---|---|---|---|
| 1 | 6.38 | 5.71 | 10.47 | 6.0099 | 5.80 | 6.3021 | 1.22 |
| 2 | 16.10 | 15.29 | 5.01 | 16.0666 | 0.20 | 15.9326 | 1.04 |
| 3 | 33.13 | 32.53 | 1.82 | 32.7091 | 1.27 | 32.3035 | 2.49 |
| 4 | 33.53 | 34.95 | 4.23 | 34.3994 | 2.59 | 34.0557 | 1.57 |
| 5 | 35.65 | 35.65 | 0.012 | 37.0980 | 4.06 | 35.8587 | 0.58 |
| 6 | 48.38 | 45.14 | 6.69 | 46.8730 | 3.11 | 48.5070 | 0.26 |
| 7 | 49.43 | 54.69 | 10.65 | 52.9663 | 7.15 | 49.5125 | 0.17 |
| 8 | 55.08 | 55.60 | 0.94 | 54.8028 | 0.50 | 54.2495 | 1.51 |
| 9 | 63.04 | 60.15 | 4.59 | 62.0316 | 1.59 | 63.5583 | 0.82 |
| 10 | 66.52 | 67.56 | 1.57 | 67.4895 | 1.45 | 67.4334 | 1.37 |
| Total Mean Error (TME) | ——— | ——— | 4.6 | ——— | 2.7 | ——— | 1.1 |

## 7. Conclusion

The TME shows significantly more accurate results obtained by the MCDWIS algorithm compared to the initial frequency and the updated frequency from the AMH simulation. The consistency of the MCDWIS algorithm with higher dimensional problems and more complex models is a testament to its validity in the practice of indirect FE model updating. It is important to note that, in the MCDWIS the weighting process and importance sampling are in control. Thus, as in literature, these are presented in the plots. Additionally, the algorithm dependent parameters are relatively simple to tune by trial and error and could be experimented with further, to achieve even more accurate results. The degrees of freedom of the element matrices, determines the dimensionality of the search space. The increase of dimensions of the system becomes relevant in our investigation due to the resultant increase of variable parameters in the system, and in turn this leads to Monte Carlo integration becoming computationally complex. Compared to other statistical methods, the MCDWIS algorithm is complex however, it has been demonstrated that it works with high accuracy and given the missing data issue when dealing with multimodal and complex structures this algorithm is useful due to the fact it relieves the need for exact or perfect sampling. The dynamic weighting and population control strategies aim to increase the convergence rate to the target distribution, through the numerical approximation of unbiased importance sampling estimates.